\documentclass[11pt]{article}
\usepackage{amsthm,amsfonts,amsmath,amssymb,graphicx}
\input xypic

\newcommand{\CC}{{\mathbb{C}}}

\newcommand{\PP}{{\mathbb{P}}}
\newcommand{\QQ}{{\mathbb{Q}}}
\newcommand{\RR}{{\mathbb{R}}}
\newcommand{\ZZ}{{\mathbb{Z}}}

\newcommand{\calO}{{\mathcal O}}

\newcommand{\calI}{{\mathcal I}}
\newcommand{\calA}{{\mathcal A}}

\newcommand{\calP}{{\mathcal P}}

\newcommand{\calX}{{\mathcal X}}

\newcommand{\op}{\operatorname}
\newcommand{\Sat}{{\calA_g^{\op {Sat}}}}
\newcommand{\Vor}{{\calA_g^{\op {Vor}}}}
\newcommand{\Perf}{{\calA_g^{\op {Perf}}}}
\newcommand{\Perfl}{{\calA_g^{\op {Perf}}(\ell)}}
\newcommand{\Part}{{\calA_g^{\op {part}}}}

\newcommand{\Prt}{{\calA_{g-1}^{\op {part}}}}
\newcommand{\X}{{\calX_{g-1}}}
\newcommand{\A}{{\calA_{g-1}}}
\newcommand{\Xl}{{\calX_{g-1}(\ell)}}
\newcommand{\Al}{{\calA_{g-1}(\ell)}}
\newcommand{\Xpart}{{\calX_{g-1}^{\op {part}}}}
\newcommand{\Xpartl}{{\calX_{g-1}^{\op {part}}(\ell)}}
\newcommand{\XPl}{{\calX_{g-1}^{\op {Perf}}(\ell)}}

\newcommand{\Centr}{{\calA_g^{\op {Centr}}}}
\newcommand{\Igu}{{\calA_g^{\op {Igu}}}}
\newcommand{\Sp}{\op{Sp}}
\newcommand{\PSp}{\op{PSp}}
\newcommand{\GL}{\op{GL}}
\newcommand{\Sym}{\op{Sym}}
\newcommand{\symm}{\op{symm}}
\newcommand{\CH}{\op{CH}}

\newcommand{\NSQ}{\op{NS}_{\mathbb Q}}

\newcommand{\Kthree}{\mathop{\mathrm {K3}}\nolimits}

\def\T{\Theta}
\def\Ti{\Theta^{(i)}}

\def\b{\beta}

\def\Abar{{\overline{\calA_g}}}

\theoremstyle{plain}
\newtheorem{thm}{Theorem}[section]
\newtheorem{lm}[thm]{Lemma}
\newtheorem{prop}[thm]{Proposition}
\newtheorem{cor}[thm]{Corollary}

\newtheorem{conj}[thm]{Conjecture}

\theoremstyle{definition}

\newtheorem{rem}[thm]{Remark}

\begin{document}
\title{Some intersection numbers of divisors on toroidal
compactifications of $\calA_g$}
\author{C. Erdenberger, S. Grushevsky\footnote{Research is
supported in part by National Science Foundation under the grant
DMS-05-55867}, and K. Hulek\footnote{Research is
supported in part by DFG grant Hu 337/5-3}} \maketitle
\begin{abstract}
We study the top intersection numbers of the boundary and Hodge
class divisors on toroidal compactifications of the moduli space $\calA_g$ of
principally polarized abelian varieties and compute those numbers that live
away from the stratum which lies over the closure of $\calA_{g-3}$ in the
Satake compactification.
\end{abstract}

\section{Introduction}
\label{sec:intro}
The moduli space $\calA_g$ of principally polarized abelian varieties
of genus $g$ has been an object of investigation for many years. It
is well known that $\calA_g$ admits various compactifications,
notably the Satake, sometimes also called the minimal
compactification $\Sat$, as well as various toroidal
compactifications. The latter depend on the choice of an admissible
fan such as the second Voronoi decomposition, the perfect cone or
first Voronoi decomposition, or the central cone decomposition. Due
to the recent work of several authors, most notably V.~Alexeev and
N.~Shepherd-Barron, it has now become more clear which role the
various compactifications play. Alexeev \cite{al} has shown that the
Voronoi compactification $\Vor$, given by the second Voronoi
decomposition, solves a moduli problem. More precisely, he has shown
that it is an irreducible component of a natural modular
compactification, and the situation has been further  clarified by
Olsson \cite{ol}. Shepherd-Barron \cite{sb} proved that the perfect
cone compactification $\Perf$ is a canonical model in the sense of
Mori theory, provided $g \ge 12$. Finally, the central cone
decomposition $\Centr$ coincides with Igusa's partial
desingularization $\Igu$ of the Satake compactification $\Sat$. The
various toroidal compactifications coincide for $g \leq 3$, but are
all different in general.

The perfect cone compactification has a very simple Picard group: it
is generated (over $\QQ$) by the determinant of the Hodge bundle $L$
and the (irreducible) boundary divisor $D$. In this paper we need to
be careful to distinguish the computations on the stack and on the
coarse moduli space: a generic ppav has an involution $-1$, and thus
a generic point of $\calA_g$ has a stabilizer of order two. Thus let
us specify explicitly that by $D$ we mean the substack
$\Perf\setminus\calA_g$ of $\Perf$.

Shepherd-Barron \cite[p. 41]{sb} has posed the question of
determining the top intersection numbers of divisors
$$
  a_N^{(g)}:=\langle L^{G-N} D^N\rangle_\Perf \eqno(*)
$$
where $G:=\frac{g(g+1)}{2}$ is the dimension of $\Perf$. It is this
question which we want to address in this paper and to which we give
a partial answer:
\begin{thm}
The only three intersection numbers with $N<3g-3$ that are non-zero
are those for $N=0,g,2g-1$ (and thus the power of $L$ being equal to
$\dim\calA_g,\dim\calA_{g-1}$, and $\dim\calA_{g-2}$, respectively).
The numbers are
\begin{equation}\label{int0}
  \langle L^{\frac{g(g+1)}{2}}\rangle_\Perf=
  (-1)^G2^{-g}G!\prod\limits_{k=1}^g\frac{\zeta(1-2k)}{(2k-1)!!}
\end{equation}
\begin{equation}\label{int1}
  \langle L^\frac{(g-1)g}{2}D^g\rangle_\Perf=\frac12
  (-1)^{G-1}(g-1)!(G-g)!\prod\limits_{k=1}^{g-1}\frac{\zeta(1-2k)}{(2k-1)!!}
\end{equation}
and
\begin{equation}\label{int2}
  \langle L^\frac{(g-2)(g-1)}{2}D^{2g-1}\rangle_\Perf=(I) + (II) + (III)
\end{equation}
where the terms (I), (II) and (III) are given explicitly by Theorems
\ref{theo:formulaforI}, \ref{theo:formulaforII} and
\ref{theo:formulaforIII} respectively.
\end{thm}
The most striking result of our computations is that
\begin{equation} \label{qu:zeroes}
  \langle L^{G-N}D^N \rangle_{\Perf}= 0 \quad {\rm unless}
  \quad G-N = \dim \calA_k \quad {\rm for\ some\ }  k\le g
\end{equation}
{\it in the range $N < 3g-3$.} This leads one naturally to
\begin{conj}
The intersection numbers $a_N^{(g)}$ {\it for any $N$} vanish unless
$G-N=k(k+1)/2$ for some $k\le g$.
\end{conj}
One can also ask this question for other toroidal compactifications
of $\calA_g$, and it is tempting to conjecture that if one interprets
$D$ as the closure of the boundary of the partial compactification,
this would also hold.

Of course, one could even hope that such a vanishing result holds for
perfect cone compactifications or in fact for all (reasonable)
toroidal compactifications of any quotient of a homogeneous domain by
an arithmetic group. This is e.g. the case for the moduli space of
polarized $\Kthree$ surfaces. However, in this case the minimal
compactification has only two boundary strata, which are of dimension
$0$ and $1$ respectively. Thus the vanishing is essentially automatic
from our discussion in Section \ref{sec:intpartial}.

Of the three non-zero numbers above, the top self-intersection number
$a_0^{(g)}=L^G$ of the Hodge line bundle $L$ can be computed using
Hirzebruch-Mumford proportionality, see e.g. \cite[Theorem
3.2]{vdg1}. The second non-zero intersection number can be easily
obtained following the methods of van der Geer from \cite{vdg2}, in
working on the partial compactification --- we give the details in
the next section. Thus in the theorem above, apart from the vanishing
results, the significant new non-zero number we compute is
$a_{2g-1}^{(g)}=\langle L^\frac{(g-2)(g-1)}{2}D^{2g-1}
\rangle_\Perf$.

The intersection theory on $\Perf$ for $g\le 3$ is described
completely (including the intersections of higher-dimensional
classes) by van der Geer \cite{vdg2}, while the intersection theory
of divisors for the perfect cone compactification as well as for the
more complicated second Voronoi compactification for $g=4$ was
computed in our previous work \cite{egh}. We have greatly benefited
from van der Geer's methods in \cite{vdg1}, \cite{vdg2}.

Throughout the paper we work with $\Perf$, but our computations also
provide information for other toroidal compactifications (see Section
\ref{sec:comments}).

\smallskip
Throughout this paper we work over the field of complex numbers
$\CC$.

\smallskip
{\footnotesize The second author would like to thank Leibniz
Universit\"at Hannover for hospitality in January 2005 and September
2006, when research for this work was performed. Both second and
third authors are grateful to the Mittag-Leffler Institute where we
had a chance to continue our collaboration in 2007.}

\section{Intersections on the partial compactification}
\label{sec:intpartial}
We shall now compute the intersection numbers
$a_N^{(g)}$ for $N \le 2g-2$. The considerations of this section are
true for all toroidal compactifications $\Abar$. We have already
remarked that the number $L^G=a_0^{(g)}$ can be computed via
Hirzebruch-Mumford proportionality. All other intersection numbers
involve the boundary and hence it is essential to obtain a good
understanding of its geometric structure. Every toroidal
compactification $\Abar$ admits a projection
\begin{equation} \label{eq:sat}
  \pi:\Abar\to\Sat =\calA_g\sqcup\calA_{g-1} \sqcup\ldots \sqcup\calA_0
\end{equation}
to the Satake compactification.

Following \cite{vdg2} we denote
$\b_k:=\pi^{-1}(\calA_{g-k}^{\op{Sat}})$. We also denote $\Part:=
\Abar\setminus\b_2$. This is the partial compactification considered
by Mumford \cite{mu}, and it is independent of the choice of the
toroidal compactification. By abuse of notation we shall also use the
symbol $D$ instead of $\b_1 \setminus \b_2$ if it is clear that we
are working away from $\b_2$. The boundary of the partial
compactification is the universal Kummer family over $\calA_{g-1}$.
More precisely, let $\calX_{g-1} \to \calA_{g-1}$ be the universal
abelian variety of dimension $g-1$ (which exists as a stack), then
there exists a degree $2$ map $j: \calX_{g-1} \to \calA_{g-1}$ and
thus (cf. \cite[p. 5]{vdg1})
\begin{equation}\label{equ:boundary}
j_*([\calX_{g-1}])=2D.
\end{equation}
Note that $D$ is denoted by $\sigma_1$ in \cite{vdg1}. In computing
the intersection numbers we have to be careful to distinguish between
stacks and varieties. We will, however, not need the full force of
the stacky machinery as all our calculations can be done by going to
suitable level covers.

In order to compute intersection numbers involving the boundary,
first note that $D\in \CH^*(\b_1)$. Here $\CH^*$ denotes the Chow
ring (tensored with $\QQ$, since torsion is irrelevant for the
intersection numbers). Hence any intersection numbers involving $D$
can be computed on $\b_1$.

The following lemma is well-known, but we include the proof for easy
reference.
\begin{lm}
For every toroidal compactification $\Abar$ we have
\begin{equation}\label{eqn:intonb1}
  L^M|_{\b_k} =0\in \CH^*(\b_k)\quad{\rm for}\quad
  M>\frac{(g-k)(g-k+1)}{2}.
\end{equation}
In particular
$$
 a_N^{(g)}=0\qquad{\rm for}\ 1\le N\le g-1.
$$
\end{lm}
\begin{proof}
Recall that suitable powers of $L$ are globally generated and can be
used to define the Satake compactification $\Sat$ as the closure of
the image of $\calA_g$ under the morphism defined by $L^{\otimes k}$.
Taking general hyperplane sections we can thus represent (a multiple
of) $L^N$ in $\Sat$ by a cycle which has empty intersection with
$\calA_{g-k}^{\op{Sat}}$, as this space has dimension
$(g-k)(g-k+1)/2$. This gives the first claim (cf. also the proof of
\cite[Lemma 1.2]{vdg1}).

The second claim follows since $D\in \CH^*(\b_1)$ and
$L^{G-N}|_{\b_1}=0$ for $N > (g-1)g/2$.
\end{proof}

In the range $N \geq 2g-2$ we know that $L^{G-N}|_{\b_2}=0$. By the
inclusion Chow exact sequence
$$
  \CH^*(\b_2)\to \CH^*(\b_1)\to \CH^*(\b_1\setminus\b_2)\to 0
$$
the intersection numbers $a_N^{(g)}$ for $g \leq N \leq 2g-2$ can be
computed on the stratum $\b_1\setminus\b_2$, i.e. on the universal
Kummer family.

It is well known (cf. \cite[Proposition 1.8]{mu} and \cite[Lemma
1.1]{vdg1}) that the pullback of the boundary divisor $D$ under the
map $j: \calX_{g-1} \to \calA_{g-1}$ is given by
\begin{equation}\label{equ:pullbackboundary}
j^*D = -2\T
\end{equation}
where $\T$ denotes the symmetric theta divisor on $\calX_{g-1}$
trivialized along the $0$-section.
Taking the
symmetric theta divisor means taking the divisor of the universal
theta function and dividing it by the theta constant. Since the theta
constant is a modular form of weight $1/2$, i.e. a section of $L/2$,
this means that $\T= \T' - L/2$ where $\T'$ is the divisor of the
universal theta function used by Mumford in \cite{mu}.

Following Mumford's ideas from \cite{mu}, it was shown by van der
Geer \cite[Formula (3)]{vdg2} that for the universal family
$\pi:\X\to\A$ the pushforwards of the powers of $\T$ are
\begin{equation}\label{push}
  \pi_*(\T^{g-1})=(g-1)![\A]\quad\text{and}\quad
  \pi_*(\T^{g+\varepsilon})=0\quad\forall \varepsilon\ge 0.
\end{equation}

{}From this one easily obtains the following intersection numbers:
\begin{prop}\label{agg}
We have
\begin{itemize}
\item[{\rm (i)}]
    $a_g^{(g)}=\frac12(-2)^{g-1}(g-1)!\,a_0^{(g-1)}=\frac12
    (-1)^{G-1}(g-1)!(G-g)!\prod\limits_{k=1}^{g-1}
    \frac{\zeta(1-2k)}{(2k-1)!!}$
\item[{\rm (ii)}]$a_{g+1}^{(g)}=\ldots=a_{2g-2}^{(g)}=0$.
\end{itemize}
\end{prop}
\begin{proof}
We have already observed that $L^{G-N}|_{\b_2}=0\in \CH^*(\b_2)$ in
the range  $g\le N<2g-1$. Using the Chow inclusion exact sequence we
obtain
$$
 a_N^{(g)}=\langle L^{G-N}D^N\rangle_\Abar
 = \frac{1}{2}\langle L^{G-N}j^* (D^{N-1})\rangle_{{\calX}_{g-1}}
$$
$$
  =\frac12\langle L^{G-N}\pi_*((-2\T)^{N-1})\rangle_\A
$$
where the factor $1/2$ in the first equality comes from equation
(\ref{equ:boundary}) and where we have used equation
(\ref{equ:pullbackboundary}) and the projection formula in the last
equality. However, we know from (\ref{push}) that
$\pi_*(\T^{g+\varepsilon})=0\in\CH^*(\A)$ for all $\varepsilon\ge 0$.
Since the intersection numbers in question can be computed on $\b_2
\setminus \b_1$, this gives claim (ii).

To obtain (i) we compute
$$
 a_g^{(g)}=\frac12(-2)^{g-1}\langle L^{\frac{(g-1)g}{2}}\pi_*(\T^{g-1})\rangle_\A
$$
$$
 =\frac12(-2)^{g-1}(g-1)!\langle L^{\frac{(g-1)g}{2}}\rangle_\A
 =\frac12(-2)^{g-1}(g-1)!a_0^{(g-1)}
$$
and thus claim (i) follows from Hirzebruch-Mumford proportionality
for genus $g-1$.
\end{proof}
\begin{rem}
These intersection numbers have been computed in a different way for
$g\le 4$ in \cite{vdg2} and \cite{egh}, and our results agree with
those --- the extra factors of $2$ in comparison with \cite{egh} are
due to the fact that we work with the stack $\calA_g$ as compared to
the variety as we did in \cite{egh}.

For example, for $g=4$ we get
$$
  -\frac{1}{7560}=\langle L^6D^4\rangle_{\overline{\calA_4}}=\frac12(-2)^3
  3! \langle L^6\rangle_{\calA_3}=-24\cdot\frac{1}{181440}.
$$
\end{rem}

\section{Geometry of the locus of corank $\le 2$ degenerations}
\label{sec:geomb3}
{}From now on we shall work with the perfect cone compactification
$\Perf$. We have already remarked that there is a map $j: \calX_{g-1}
\to \Part$ which maps the universal family $2:1$ onto the boundary.
It is crucial for us that this picture of the boundary can be
extended one step further. More precisely, the universal family $j:
\calX_{g-1} \to \calA_{g-1}$ has a partial compactification $\Xpart
\to \Prt$ such that the map  $j: \calX_{g-1} \to \Part$ extends to a
map $j: \Xpart \to \Perf$ with image $j_*[\Xpart]=2[\b_1 \setminus
\b_3]$.

For details of the geometry (of the underlying variety) of the stack
$\Xpart$ see \cite{ts}, \cite{hu}. Here we recall the essential
properties of this stack. The fibration $\pi: \Xpart \to  \Prt$
extends the universal family of genus $g-1$. The fiber over a point
in  $\calX_{g-2} = \Prt \setminus  \A$ is a corank $1$ degeneration
of a $(g-1)$-dimensional abelian variety, i.e. a $\PP^1$-bundle over
an abelian variety of genus $g-2$ where the $0$-section and the
$\infty$-section are glued with a shift. We denote the composition of
the projection $\Xpart \to \calX_{g-2}$  with the projection
$\calX_{g-2} \to \calA_{g-2}$ by $\pi_{\op{Sat}}$. Then we have the
following picture: let $B\in\calA_{g-2}$ and let $\calP_B$ denote the
Poincar\'e line bundle on $B\times B$. We denote the coordinates on
$B \times B$ by $(z,b)$. Note that the coordinates are {\it not}
interpreted symmetrically --- the first one defines an actual point
on $B$, while the second is thought of as parameterizing the moduli
of semiabelian varieties with abelian part $B$. Then
\begin{equation}\label{bundle}
  \pi_{\op{Sat}}^{-1}(B)=\PP^1(\calP_B\oplus\calO)/(z,b,0)\sim (z+b,b,\infty).
\end{equation}
For a discussion of this see also \cite[p. 356]{mu}. We would like to
point out, however, that in our convention the role of the two
factors is interchanged compared to \cite{mu}.

The map  $\pi: \Xpart \to  \Prt$ is no longer smooth, but drops rank
at the singular locus of the degenerate abelian varieties. We denote
this locus, which will be crucial to our considerations, by $\Delta$.
Note that $\Delta = \calX_{g-2}\times_{\calA_{g-2}}\hat{\calX}_{g-2}$
(see e.g.~\cite[p.~11]{vdg2}). Using the principal polarization we
can identify this with
${\calX_{g-2}}\times_{\calA_{g-2}}{\calX_{g-2}}$, but we will be
careful to usually not perform this identification, as the two
factors in $\Delta$ geometrically play a very different role. Clearly
$\Delta$ is of codimension two in $\Perf$. As a stack it is smooth.

We are also interested in the restriction of the above family to the
boundary, i.e. in the family $\pi: \Xpart \setminus \calX_{g-1} \to
\calX_{g-2}$. The total space is now no longer smooth, in fact it is
a non-normal with singular locus $\Delta$. Locally along
$\Delta$ we have two smooth divisors intersecting transversally. If
we remove $\Delta$ from $\Xpart$ we obtain a $\CC^*$-fibration over
$\calX_{g-2}\times_{\calA_{g-2}}\hat{\calX}_{g-2}$ whose total space
can be identified with the Poincar\'e bundle $\calP$ with the
$0$-section removed. If $Y$ is the normalization of $\Xpart$, then
this implies that $Y=\PP^1(\calP\oplus\calO)$. We must keep track of
the various projection maps involved and therefore want to give them
names. The geometry is summarized by the following commutative
diagram of projection maps:
\begin{equation}\label{project}
  \xymatrix{Y=&\PP^1(\calP\oplus\calO)\ar^f[d]\ar_{\tilde\pi}[ddl]& \\
  \Delta=&\qquad\calX_{g-2}\times_{\calA_{g-2}}\hat{\calX}_{g-2}
  \ar^{pr_1}[dl]\ar_{pr_2}[dr]\ar^h[dd]& \\
  \calX_{g-2}\ar^{p_2}[dr]&& \hat{\calX}_{g-2}\ar_{p_1}[dl]\\
  &\calA_{g-2}& }
\end{equation}
Here the notation for the map $\tilde\pi$ is due to the fact that it
is the composition of the normalization map $Y\to\Xpart$ with the
restriction of the universal family ${\pi}:\Xpart \to \Prt$ to the
boundary $\calX_{g-2}$ of $\Prt$.

\section{Intersection theory on $\Delta$ and $Y$}
\label{sec:inttheorydeltaY}
In the beginning of this section the geometry of $\Delta$ sitting
inside $\Perf$ will not be important for us, and thus we will for now
use the principal polarization to identify $\Delta=\calX_{g-2}
\times_{\calA_{g-2}} \calX_{g-2}$. We first want to understand the
N\'eron-Severi group of $\QQ$-divisors on $\Delta$ modulo numerical
equivalence. Let $T_i=pr_1^*(\T)$ for $i=1,2$ where $pr_i$ is the
$i$-th projection from $\calX_{g-2}\times_{\calA_{g-2}} \calX_{g-2}$
to $\calX_{g-2}$ {}(as in the diagram above), and let $\calP$ be the
Poincar\'e bundle on $\calX_{g-2} \times_{\calA_{g-2}} \calX_{g-2}$.

\begin{lm} \label{lm: NeronSeverigroup}
For $g \geq 4$ the N\'eron-Severi group  of $\QQ$-divisors on
$\Delta$ modulo numerical equivalence has rank $4$; more precisely
$$
  \NSQ(\Delta)=\NSQ(\calX_{g-2} \times_{\calA_{g-2}} \calX_{g-2})=
  \QQ h^*L\oplus \QQ T_1\oplus\QQ T_2\oplus\QQ\calP.
$$
If $g=3$ then $\NSQ(\calX_{g-2} \times_{\calA_{g-2}}\calX_{g-2})$ has
rank $3$ and is generated by $T_1, T_2$ and $\calP$.
\end{lm}
\begin{proof}
We first remark that for the basis of the fibration $\calX_{g-2}
\times_{\calA_{g-2}}\calX_{g-2} \to \calA_{g-2}$ we have
$\NSQ(\calA_{{g-2}})=\QQ L$, provided $g \geq 4$, otherwise it is
trivial. We also observe that for a very general fiber $A \times A$
(i.e. outside a countable union of proper subvarieties) the rank of
$\NSQ(A \times A)$ equals $3$ \cite[Chapter 5.2]{bila}. The fibration
$\calX_{g-2} \times_{\calA_{g-2}} \calX_{g-2}\to \calA_{g-2}$ is
topologically locally trivial. If $M$ is a line bundle on
$\calX_{g-2} \times_{\calA_{g-2}} \calX_{g-2}$ then we can write
$$
  c_1(M)= c_0h^*L + c_1[T_1] + c_2[T_2] + c_3[\calP]+
  \sum_{i=4}^{2\binom{2g-4}{2}}c_i[\lambda_i]
$$
where the $c_i$ are rational numbers and the $\lambda_i \in \ H^2(A
\times A, \QQ)$ are chosen such that they form a basis of $H^2(A
\times A, \QQ)$ together with  $[T_1], [T_2]$ and $[\calP]$.

The coefficients of this sum are locally constant. Restricting to a
very general fiber, we see that $c_i=0$ for $i \geq 4$.

For an alternative proof see \cite[Lemma 2.10]{sb}.
\end{proof}

We will now want to understand the intersection theory of divisors on
$\calX_{g-2}\times_{\calA_{g-2}}\calX_{g-2}$. Much of what we need
was done in \cite{gl}, where one fiber of this universal family was
studied. Since we need some more information and our notations are
different, we review and extend the setup.

We have a shift operator
$$
  s:  \calX_{g-2} \times_{\calA_{g-2}} \calX_{g-2} \to  \calX_{g-2}
  \times_{\calA_{g-2}} \calX_{g-2}
$$
$$
 s(\tau,z,b):= (\tau,z+b,b)
$$
where we will use notations $z$ and $b$ for the variables on the
first and second factor of $\calX_{g-2} \times_{\calA_{g-2}}
\calX_{g-2}$, respectively. It is useful to keep in mind that
geometrically $b$ lies on the dual abelian variety.

\begin{lm}[\cite{gl}] \label{lm:shiftoperator}
The shift operator acts on $\NSQ(\calX_{g-2} \times_{\calA_{g-2}}
\calX_{g-2})$ by
\begin{itemize}
\item[{\rm (i)}] $s^*T_1=T_1 + T_2 + \calP$
\item[{\rm (ii)}] $s^*T_2=T_2$
\item[{\rm (iii)}] $s^*\calP = 2T_2 + \calP$
\item[{\rm (iv)}] $s^*(h^*L)=h^*L$.
\end{itemize}
\end{lm}
\begin{proof}
Equality (iv) is obvious, since $h\circ s=h$. We shall prove this
lemma by intersecting with suitable test curves. For this let
$A\in\calA_{g-2}$ be a general ppav, and let $C\subset A$ be a
general curve of degree $n:=\langle C\T_A\rangle_A$ (where $\T_A$
denotes the theta divisor on $A$); denote by $i_1,i_2,$ and $i_d$ the
three inclusions $A\hookrightarrow A\times A$ given by mapping to the
first factor, to the second factor, and to the diagonal respectively.
Moreover let $C_L\subset\calA_{g-2}$ be a generic curve, which we
think of as sitting in the $0$-section of $\calX_{g-2}
\times_{\calA_{g-2}} \calX_{g-2} \to \calA_{g-2}$. The intersection
matrix of these curves with the generators of $\NSQ(\calX_{g-2}
\times_{\calA_{g-2}} \calX_{g-2})$ is given by
\begin{equation}
 \begin{array}{|c|cccc|}
 \hline
 {}_{\rm divisor}\!\!\backslash^{\rm curve}&C_L& C_1&C_2&C_d\\ \hline
 h^*L    &\langle C_L L\rangle_{\calA_{g-2}}&0&0&0\\
 T_1  &0&n&0&n\\
 T_2  &0&0&n&n\\
 \calP&0&0&0&2n\\
 \hline
 \end{array}
\end{equation}
The intersection of $h^*L$ with the curves $C_1$, $C_2$ and $C_d$ is
clearly $0$ since these curves lie in a fiber of the map
$h:\calX_{g-2} \times_{\calA_{g-2}} \calX_{g-2} \to \calA_{g-2}$ .
The intersections of $T_i$ with $C_1,C_2,$ and $C_d$ are obvious. The
intersection of $T_i$ and $\calP$ with $C_L$ is $0$ since these
bundles are, by definition, trivial along the $0$-section which
contains $C_L$. Lastly, we note that the Poincar\'e bundle is trivial
on $A\times \lbrace 0\rbrace$ and $\lbrace 0\rbrace\times A$, and
thus trivial on $C_1$ and $C_2$; however, on the diagonal the
Poincar\'e bundle restricts to $\calO_A(2\T_A)$ (see \cite{mu}, p.
357). {}The intersection {}matrix above is non-singular; hence the
curves $C_1$, $C_2$, $C_d$ and $C_L$ can serve as test curves for the
intersection theory on $\calX_{g-2} \times_{\calA_{g-2}}
\calX_{g-2}$.

We shall now calculate the intersection of $(s^N)^*(T_1)$ with the
test curves where $N \geq 1$ is an integer. Recall that $T_1$ is
given by the pullback of the divisor $\lbrace
\Theta(\tau,z)/\Theta(\tau,0)=0 \rbrace$ on $\calX_{g-2}$ via $pr_1$.
Here $\Theta(\tau,z)=\Theta_{00}(\tau,z)$ is the standard theta
function and the denominator $\Theta(\tau,0)$ ensures that $T_1$ is
trivial along the $0$-section. Applying the shift operator $N$ times
means that we have to consider the zeroes of the function
$\Theta(\tau,z+Nb)/\Theta(\tau,Nb)=0$. Note that this is still
trivial along the $0$-section, hence the intersection with $C_L$
remains trivial. To compute the intersection with the curves $C_i$
for $i=1,2,d$ the denominator is irrelevant as these curves are
contained in a fixed fiber. We claim that the number of zeroes of
$\Theta(\tau,z+Nb)$ on $C_1$, $C_2$ and $C_d$ is $n$, $N^2n$ and
$(N+1)^2n$ respectively. The first number follows since $C_1$ is
contained in $b=0$ and hence the number of zeroes does not change.
For $C_2$ we have to compute the number of zeroes of
$\Theta(\tau,Nb)$ where the variable is now $b$. The claim follows
since multiplication by $N$ induces multiplication by $N^2$ on the
second cohomology of an abelian variety, and hence the intersection
number is multiplied by $N^2$. For $C_d$ we notice that we have to
put $z=b$ and the claim for this curve now follows in the same way.
{}From these intersection numbers we deduce that $(s^N)^*(T_1)=T_1 +
N^2T_2 + N\calP$. Putting $N=1$ gives claim (i). Claim (ii) follows
since the shift operator $s$ leaves the second variable unchanged.

To prove (iii) we can use the formulae $(s^2)^*(T_1)=T_1 + 4T_2 +
2\calP$ (put $N=2$) and
$(s^2)^*(T_1)=s^*(T_1+T_2+\calP)=(T_1+T_2+\calP) + T_2 + s^*\calP$.
Comparing these two formulae gives (iii).
\end{proof}
We now return to the discussion of the geometry of $\Xpart$ and its
normalization $Y=\PP(\calP \oplus \calO)$. We shall need to know the
normal bundle of $\Delta$ in $\Xpart$.
\begin{prop}\label{normaltoDelta}
The normal bundle
$$
N_{\Delta/\Xpart}=\calP\oplus (\calP^{-1}\otimes T_2^{-2})
$$
i.e. $\Delta^2=\Delta|_\Delta=-\calP\cdot(\calP+2T_2)|_\Delta$
in $\CH^*(\Xpart)$.
\end{prop}
\begin{proof}
For smooth embeddings $X\subset Y\subset Z$ we have
$$
  0\to N_{X/Y}\to N_{X/Z}\to N_{Y/Z}|_X\to 0
$$
and thus the normal bundle of an intersection of two varieties is the
sum of their normal bundles. In our situation the total space
$\Xpart$ is smooth, containing $\Delta$, which is also smooth (all of
this in a stack sense). Locally, $\Delta$ is a local complete
intersection of the two branches of $Y$ glued along the $0$-section
and the $\infty$-section. Thus the normal bundle of $\Delta$ is the
sum of $N_{\rm 0-section/Y}$ and of the pullback of $N_{\infty-{\rm
section/Y}}$ under the shift which controls the gluing of the
$0$-section and the $\infty$-section. The normal bundle to the
$0$-section is just the bundle we have over it, i.e. the Poincar\'e
bundle. This gives the first summand. We can reverse the role of the
$0$-section and the $\infty$-section by considering $\PP(\calO \oplus
\calP^{-1})$ instead of $\PP(\calP \oplus \calO)$ which gives the
second summand. The normal bundle is thus a direct sum of $\calP$ and
$s^*(\calP^{-1})=\calP^{-1} \otimes T_2^{-2}$ where the last equality
now follows from Lemma \ref{lm:shiftoperator}.

The second claim then follows immediately from the double point
formula.
\end{proof}

For what follows we need to understand the intersection theory of the
$\PP^1$-bundle $Y$. We will denote the class of the $0$-section by {}
$\xi\in\NSQ(Y)$. Then the class of the $\infty$-section is {} $\xi -
f^*\calP\in\NSQ(Y)$.

\begin{lm}
The group of divisors on $Y$ up to numerical equivalence is
$\NSQ(Y)=f^*\NSQ(\Delta)\oplus \QQ\xi$. Moreover
\begin{equation}\label{xisquare}
  \xi(\xi-f^*\calP)=0
\end{equation} in $\CH^*(Y)$.
\end{lm}
\begin{proof}
Clearly  $\NSQ(Y)$ is generated by $f^*\NSQ(\Delta)$ together
with the class $\xi$. The relation $\xi(\xi-f^*\calP)=0$ is then
simply the fact that the $0$-section and the $\infty$-section do not
intersect. See also \cite[Remark 3.24]{fu}.
\end{proof}

\begin{prop}[see \cite{gl}] \label{prop:intersectionY}
The pullback of the symmetric theta divisor $\T$ to $Y$ is numerically
\begin{equation}\label{thetaonY}
 \T|_Y=\xi+f^*T_1-\frac12 f^*\calP\in \NSQ(Y).
\end{equation}
\end{prop}
\begin{proof}
We refer to \cite{gl} for a more general discussion of divisor
classes on the universal semiabelian family, possibly with a level
structure. Since we are interested in the universal symmetric theta
divisor (as opposed to the universal theta function), for
completeness we give a full proof below. We prove this result again
by computing the intersections of the divisor classes with suitable
test curves.

We use the curves $C_1$, $C_2$, $C_d$ and $C_L$ as in the proof of
Lemma \ref{lm:shiftoperator}, but this time we embed them into the
$\infty$-section of $Y$. As a fifth test curve we consider a general
fiber $F$ of the $\PP^1$-bundle $Y$. We claim that the $5$ curves
described are a suitable set of test curves. Indeed, this follows
from the following intersection matrix
\begin{equation}
 \begin{array}{|c|ccccc|}
 \hline
 {}_{\rm divisor}\!\!\backslash^{\rm curve}&C_L& C_1&C_2&C_d&F\\ \hline
 f^*h^*L    &\langle C_L L\rangle_{\calA_{g-2}}&0&0&0&0\\
 f^*T_1  &0&n&0&n&0\\
 f^*T_2  &0&0&n&n&0\\
 f^*\calP&0&0&0&2n&0\\
 \xi     &0&0&0&0&1\\
 \hline
 \end{array}
\end{equation}
The top left hand $4 \times 4$ matrix follows immediately from Lemma
\ref{lm:shiftoperator}. The intersection of the divisors $f^*h^*L$,
$f^*T_1$, $f^*T_2$ and $f^*\calP$ with $F$ is clearly $0$ since these
curves are contracted by $f$. Similar the intersection of $\xi$ with
these curves is $0$ since $\xi$ is the $0$-section and these curves
lie in the $\infty$-section. The remaining number $\xi.F=1$ is also
obvious.

We now want to compute the intersection numbers of the symmetric
theta divisor with these test curves. The symmetric theta divisor on
the universal semiabelian variety is given by the function
$$
  \Theta_{\symm}=\frac{\theta(\tau,z+ b/2) + x\theta(\tau,z - b/2)}{2\theta(\tau,b/2)}
$$
where $\theta(\tau,z)=\Theta_{00}(\tau,z)$ is the standard theta
function. Here $x$ is the variable on the $\PP^1$ fiber of $Y$ over
$\Delta$. Note that the denominator ensures that the corresponding
divisor is trivial on the section given by the neutral element of the
universal semiabelian variety, since setting $z=0$ and $x=1$ gives
the constant $1$ (see also \cite{huwe}). Clearly the degree of
$\Theta_{\symm}$ on $F$ is $1$. The curves $C_1$, $C_2$, $C_d$ and
$C_L$ all lie in the $\infty$-section of $Y$. To evaluate
$\Theta_{\symm}$ on these curves therefore means that we have to
consider only the second summand ${\Theta(\tau,z -
b/2)}/{2\Theta(\tau,b/2)}$. To compute the intersection with $C_L$ we
have to set $z=0$ and we immediately find that this degree is $0$.
The same argument applies to $C_2$. For $C_1$ we must put $b=0$ and
we see that the degree equals $n$ (recall that our computations
happen on a general abelian variety, and in particular that we can
assume $\Theta(\tau,0)\ne 0$, so that the denominator is non-zero).
Finally, in order to compute the intersection number with $C_d$ we
put $z=b$ and we obtain again $0$. The discussion above is summarized
by the following intersection numbers:
\begin{equation}
 \begin{array}{|c|ccccc|}
 \hline
 {}_{\rm divisor}\!\!\backslash^{\rm curve}&C_L& C_1&C_2&C_d&F\\ \hline
 \T|_Y&0&n&0&0&1\\
 \hline
 \end{array}
\end{equation}
It now follows immediately from these intersection numbers that
$$
 \T|_Y=\xi+f^*T_1-\frac12 f^*\calP
$$
as claimed.
\end{proof}

\begin{cor}\label{thetaondelta}
We have
$$
\T|_\Delta= T_1 +\frac12 \calP\in \NSQ(\Delta)
$$
\end{cor}
\begin{proof}
This follows immediately from Proposition \ref{prop:intersectionY} by
restricting to the $0$-section since $\xi^2= \xi f^*\calP$ (Note that
we obtain the same result if we first restrict to the
$\infty$-section and then pull back via $s$.)
\end{proof}

\section{Level cover of the moduli space}
\label{sec:level}
In order to do our intersection calculation we want to pass to a
suitable level cover of $\calA_g$. We shall only consider level
covers of full level $\ell$ and we assume $\ell \geq 3$ to be prime.
This simplifies many of the formulae and is clearly sufficient. We
denote by $\Sat(\ell)$ and  $\Perfl$ the Satake and the perfect cone
compactification of $\calA_g(\ell)$ respectively. They fit into a
commutative diagram
$$
  \xymatrix{\Perfl  \ar^\sigma[r]\ar^\pi[d]&\Perf\ar^\pi[d] \\
  \Sat(\ell)\ar^\sigma[r] &\Sat}
$$
where $\sigma$ is a Galois cover with Galois group $\Sp(2g,\ZZ/\ell
\ZZ)$. The above diagram is interpreted as a diagram of stacks. Note
that for level $\ell \geq 3$ the element $-1$ is no longer contained
in the level-$\ell$ subgroup of $\Sp(2g,\ZZ)$ and we can, therefore,
compute intersection numbers on the level covers directly on the
underlying varieties.

We will denote by $\b_i(\ell)$ the preimage in $\Perfl$ of
$\calA_{g-i}^{\op {Sat}}(\ell)$. Unlike the no level case, the
boundary of $\Perfl$  has many irreducible components, which we
denote $D_i(\ell)=D_i$. We also recall (cf. \cite[p. 4]{vdg1}) that
\begin{equation} \label{equ:branching}
  \sigma^*(D)=\sum_i \ell D_i.
\end{equation}
The geometry of toroidal compactifications is controlled by the
admissible fan which is chosen. In the case of $\calA_g$ this is a
fan in the space $\Sym^{\geq 0}(g,\RR)$ of real semi-positive
definite symmetric $g \times g$ matrices. Here we work with the
perfect cone or first Voronoi decomposition, which is defined by
taking the convex hull of all primitive rank one forms in $\Sym^{\geq
0}(g,\ZZ)$. We want to understand the geometry of $\Perfl\setminus
\b_3(\ell)$. For this we only need the intersection of the perfect
cone decomposition with the subspace $\Sym^{\geq 0}(2,\RR) \subset
\Sym^{\geq 0}(g,\RR)$ where the inclusion is given by taking the
first two variables. For $g \leq 3$ the first and the second Voronoi
decomposition coincide, and both decompositions are well understood
for $g \leq 4$, see \cite{v1}, \cite{v2a} and \cite{v2b} as well as
\cite{hs} for a more recent exposition. Note that as in the case
without level structure
$$
D_i \cong \Xpartl.
$$
We first recall
\begin{prop}\label{Disthetalevel}
The following identity holds in $\CH^*(\Perfl \setminus \b_3(\ell))$:
$$
D_i|_{D_i}=-2\frac{\T}{\ell}|_{D_i}.
$$
Here $\T$ is the pullback of the theta divisor trivialized along the
$0$-section to the universal family $\Xpartl$ with level-$\ell$
structure.
\end{prop}
\begin{proof}
This follows from \cite[Proposition 3.2]{hu}.
\end{proof}

We will now discuss the intersection of the boundary components.
Recall that the boundary components  $D_i$ correspond to primitive
integer quadratic forms of rank $1$, i.e. to squares
$q_i(x)=l_i^2(x)$ where the $l_i$ are primitive integral linear
forms. The intersection of two different components $D_i$ and $D_j$
is always contained in $\beta_2(\ell)$, and two components $D_i$ and
$D_j$ have a non-empty intersection if and only if the corresponding
quadratic forms $l_i^2$ and $l_j^2$ are part of a basis of
$\Sym(g,\ZZ)$. The situation is more complicated for the
intersection of three different boundary components $D_i$, $D_j$ and
$D_k$. In this case we have two possibilities: the linear forms
$l_i$, $l_j$ and $l_k$ can either be linearly independent or not. In
the first case the intersection of the three boundary components is
contained in $\b_3(\ell)$. We will refer to this as the {\it global}
case. We will not be concerned with this as this will not contribute
to the intersection numbers. The other case is called the {\it
local} case and we will have to consider these intersections.

The intersection of four boundary divisors will never contribute to
the intersection numbers that we will compute. This follows from
\begin{lm}\label{no4intersect}
For any four distinct indices $i,j,k,l$ we have
$$
  D_i\cap D_j\cap D_k\cap D_l\subset\b_3(\ell)
$$
and thus $D_iD_jD_kD_l=0\in \CH^*(\Perf(\ell)\setminus\b_3(\ell))$.
\end{lm}
\begin{proof}
This follows from the fact that there are no cones in the perfect
cone decomposition of  $\Sym^{\geq 0}(2,\RR)$ which are spanned by
$4$ vectors.
\end{proof}
We now want to understand the intersection of divisors more
systematically. We will have to compute intersections of the form
$D_i^aD_j^bD_k^cL^d$ such that $a+b+c+d=G$. For symmetry reasons
these numbers will only depend on the integers $a$, $b$ and $c$. To
compute these intersections we will, however, work on a specific
boundary component $D_i$ and thus we must understand how $D_i|_{D_i
\cap D_j}$ and $D_j|_{D_i \cap D_j}$ compare. Boundary components
correspond to rank $1$ symmetric forms. As we are only interested in
the case where the intersection is of local type, these three forms
are linearly dependent. Modulo the action of the group $\GL(g,\ZZ)$
we can assume that the boundary components $D_i$, $D_j$ and $D_k$
correspond to the matrices
$\begin{pmatrix} 1&0\\ 0& 0\end{pmatrix}$,
$\begin{pmatrix} 0&0\\ 0& 1\end{pmatrix}$ and $\begin{pmatrix} 1&-1\\
-1& 1\end{pmatrix}$ respectively. This is shorthand for $g \times
g$-matrices where these matrices sit in the upper left corner and the
rest of the entries are zero. Correspondingly we consider
$\GL(2,\ZZ)$ as a subgroup $\GL(g,\ZZ)$. An element $\gamma \in
\GL(g,\ZZ)$ acts on a matrix $M$ by $M \mapsto
{}^t\gamma^{-1}M\gamma^{-1}$.

Now consider the element $\gamma= \begin{pmatrix} 0&1\\ 1&
0\end{pmatrix}$ interpreted as an element in $\GL(g,\ZZ)$. This
interchanges $D_i$ and $D_j$ while it fixes $D_k$. In terms of the
symplectic group we are considering the matrix
$$  \gamma=
\begin{pmatrix} 0&1 & & & & \\
 1& 0 & & & &\\
  &  &  {\bf{1}_{g-2}}  & & & \\
  &  &         & 0 & 1 & \\
  &  &         & 1 & 0 & \\
  &  &         &   &   & {\bf{1}_{g-2}}
\end{pmatrix}
$$
This matrix acts on Siegel space by
$$
\gamma
\begin{pmatrix} \tau_{11}& \tau_{12} & z\\
\tau_{21}  & \tau_{22}  & b \\
z  & b  &  \tau_{g-2}
\end{pmatrix}
=
\begin{pmatrix} \tau_{22}& \tau_{21} & b\\
\tau_{12}  & \tau_{11}  & z \\
b  & z  &  \tau_{g-2}
\end{pmatrix}.
$$

The non-empty intersection $Y_{ij}=D_i \cap D_j$ of two boundary
components is irreducible, since its toric part is a torus orbit.
Geometrically, it is a $\PP^1$-bundle over $\calX_{g-2}(\ell)
\times_{\calA_{g-2}(\ell)} \hat\calX_{g-2}(\ell)$. More precisely, it
is the compactification of the pull back $\sigma^*(\cal P)$ under the
Galois cover (this pullback is in fact the $\ell$-th power of
Poincar\'e bundle on $\calX_{g-2}(\ell) \times_{\calA_{g-2}(\ell)}
\hat\calX_{g-2}(\ell)$, see \cite{gl}, but this does not matter for
us).
\begin{lm}\label{lm:symmetricintersection} We have
\begin{itemize}
\item[{\rm (i)}] $\ell D_i|_{D_i \cap D_j} = - 2\xi - 2f^*T_1 + f^*\calP$
\item[{\rm (ii)}]  $\ell D_j|_{D_i \cap D_j} = - 2\xi - 2f^*T_2 + f^*\calP$
\end{itemize}
Here we denote a class on $Y$ and its pullback to the Galois cover
$Y_{ij}=D_i \cap D_j$ by the same symbol. The same applies to the map
$f$ which is the analogue of the bundle projection in diagram
(\ref{project}).
\end{lm}
\begin{proof}
To clarify notation we denote the trivialized theta divisor on $D_i$
by $\T^{(i)}$. It follows from Proposition \ref{prop:intersectionY}
that $\T^{(i)}|_{D_i \cap D_j}= \xi + f^*T_1 - 1/2 f^*\calP$. Here $\xi=D_i \cap
D_j \cap D_k$. We thus obtain claim (i) from Proposition
\ref{Disthetalevel} i.e.
$$
 \ell D_i|_{D_i \cap D_j} = -2\T^{(i)}|_{D_i \cap D_j} = - 2\xi - 2f^*T_1 + f^*\calP.
$$
Claim (ii) now follows from
$$
\ell D_j|_{D_i \cap D_j} =\gamma^*(\ell D_i)|_{D_i \cap D_j}
$$
$$
=(\gamma|_{D_i \cap D_j})^*(\ell D_i|_{D_i \cap D_j})
=(\gamma|_{D_i \cap D_j})^*(-2\T^{(i)}|_{D_i \cap D_j})
$$
$$
=(\gamma|_{D_i \cap D_j})^*( - 2\xi - 2f^*T_1 + f^*\calP) =  - 2\xi - 2f^*T_2 + f^*\calP
$$
where the last equality follows since $\gamma$ interchanges $z$ and
$b$.
\end{proof}

The preimage of $\Delta$ under the Galois cover is also reducible,
more precisely $\Delta(\ell)=\cup \Delta_{ijk}$ where
$\Delta_{ijk}=D_i \cap D_j \cap D_k$ with $\Delta_{ijk} \cong
\calX_{g-2}(\ell) \times_{\calA_{g-2}(\ell)} \hat\calX_{g-2}(\ell)$.
We also note
\begin{lm}\label{lm:intthreeboundcomp}
We fix an identification $D_i=\calX_{g-1}(\ell)$ and thus also an
identification $\Delta_{ijk}=D_i\cap D_j\cap D_k=
\calX_{g-2}(\ell)\times_{\calA_{g-2}(\ell)}\hat\calX_{g-2}(\ell)$.
Then the following holds:
\begin{itemize}
\item[\rm{(i)}] $D_i|_{D_i \cap D_j \cap D_k}= -2T_1 - \calP$
\item[\rm{(ii)}] $D_j|_{D_i \cap D_j \cap D_k}= -2T_2 - \calP$
\item[\rm{(iii)}]  $D_k|_{D_i \cap D_j \cap D_k}= \calP$
\end{itemize}
where we again use the same symbol for a class on $\Delta$ and its pullback under the
Galois cover.
\end{lm}
\begin{proof}
To see the first identity we have to restrict the equality $D_i|_{D_i
\cap D_j} = - 2\xi - 2f^*T_1 + f^*\calP$ to $D_i \cap D_j \cap D_k =
\xi$. The claim then follows since $\xi^2=\xi f^*\calP$. The second
identity follows in the same way from Lemma \ref{lm:symmetricintersection}.
The last claim follows since $D_k|_{D_i \cap D_j \cap D_k}= \xi|{\xi}
= \calP$.
\end{proof}
\begin{rem}
Note that (iii) is also consistent with the fact that
$N_{\Delta/Y}=N_{D_i \cap D_j \cap D_k/D_i \cap D_j}= D_k|_{D_i \cap
D_j \cap D_k}= \calP$.
\end{rem}

\section{Combinatorics of the level cover}
\label{sec:combinatoricscover}
We will now need to deal more carefully with the combinatorics of the
boundary of the level cover, to understand how the intersection
numbers that we want to compute, which are of course well-defined on
$\Perf$ with no level, can first be computed by doing an honest
computation on the variety $\Perfl$, and then how this computation on
$\Perfl$ can, due to its invariance under the deck transformation
group of the cover, be reduced back to some computations without a
level.

\smallskip
We denote by $\sigma:\Perfl\to \Perf$ the level cover, which has degree
\begin{equation}\label{nug}
  \nu_g(\ell)= |\Sp (2g,\ZZ/\ell\ZZ)|= \ell^{g(2g+1)}(1-l^{-2})\cdot\ldots\cdot(1-l^{-2g}).
\end{equation}
Note that this is the degree of the map between the stacks. The
degree of the maps between the varieties is $|\PSp (2g,\ZZ/\ell\ZZ)|$
which is $\nu_g(\ell)/2$. For further use, we denote by $d_g(\ell)$
the number of irreducible components $D_i$ of the boundary $\Perfl$.
This is equal to one half the number of points in $(\ZZ/\ell\ZZ)^{2g}
\setminus \{0 \}$, which is
\begin{equation}\label{dg}
  d_g(\ell)=\frac12\ell^{2g}(1-\ell^{-2g}).
\end{equation}
The cover $\sigma$ branches to order $\ell$ along the boundary, cf. (\ref{equ:branching}).
Using Lemma \ref{no4intersect} we get
$$
 a_N^{(g)}=\langle L^{G-N} D^N\rangle_\Perf=
 \frac{1}{\nu_g(\ell)}\langle \sigma^*L^{G-N} \sigma^*D^N\rangle_\Perfl
$$
$$
 =\frac{\ell^N}{\nu_g(\ell)}\left\langle \sigma^*L^{G-N}\left[\sum\limits_i D_i^N+
 \sum\limits_{i > j;\ a+b=N, a,b>0} \binom{N}{a}D_i^a D_j^b\right.\right.
$$
\begin{equation}\label{levelintersection}
 \left.\left.\qquad\qquad
 +\sum\limits_{i > j > k;\ a+b+c=N, a,b,c>0} \binom{N}{a,\ b,\ c}D_i^a D_j^b D_k^c
 \right]\right\rangle_{\Perfl}
\end {equation}
(where we remind the reader that for $a+b+c=N$ by definition $\binom{N}{a,\ b,\ c}
=\frac{N!}{a!b!c!}$).
We will now reduce this computation of the intersection numbers on a
level cover to the intersection computation on the no level moduli
space. Let us start with the first term.
\begin{lm} \label{lm:termIred}
The first term of equation (\ref{levelintersection}) is given by
$$
  \frac{\ell^N}{\nu_g(\ell)}\left\langle \sigma^*L^{G-N}
   \sum\limits_i D_i^N\right\rangle_\Perfl =
  \frac12\left\langle L^{G-N}(-2\T)^{N-1}\right\rangle_\Xpart.
$$
\end{lm}
\begin{proof}
We denote by $e_{I}(\ell)$ the degree of the map (as stacks) from
each boundary component $D_i$ of $\Perfl$ to the (closure of the) universal family
$\X$ (this map is clearly onto). {}From the description of $D_i$ as the
universal level family $\Xl$, mapping to $\X$, we get two factors for
this degree --- one is the degree of the map $\Al$ to $\A$ (thus it is
$\nu_{g-1}(\ell)$), and the other factor is the degree of the map of
a single fiber of the family $\Xl\to\Al$ to a fiber of $\X\to\A$, i.e
of the map of an $\ell$-times principally polarized abelian variety
of dimension $g-1$ to a principally polarized abelian variety of
dimension $g-1$. We thus have
$$
 e_{I}(\ell)=\nu_{g-1}(\ell)\ell^{2g-2}.
$$
By Proposition \ref{Disthetalevel} we obtain
$$
 \frac{\ell^N}{\nu_g(\ell)}\langle \sigma^*L^{G-N}\sum\limits_i D_i^N\rangle_\Perfl =
\frac{\ell}{\nu_g(\ell)}\langle
\sigma^*L^{G-N} \sum\limits_i (-2\Ti)^{N-1}D_i\rangle_\Perfl
$$
$$
=\frac{\ell}{\nu_g(\ell)}\, d_g(\ell)e_{I}(\ell)\langle L^{G-N} (-2\T)^{N-1} \rangle_\Xpart.
$$
To prove the lemma it remains thus to compute, substituting the
expressions for $d_g$ and $\nu_g$ from formulae (\ref{dg}) and
(\ref{nug}),
$$
 \frac{\ell}{\nu_g(\ell)}\, d_g(\ell)e_{I}(\ell)
$$
$$
=\frac{\frac12\ell\ell^{2g}(1-\ell^{-2g})\frac12\ell^{(g-1)(2g-1)}(1-\ell^{-2})\cdot
\ldots\cdot(1-\ell^{-2g+2})\ell^{2g-2}}{\frac12\ell^{g(2g+1)}(1-\ell^{-2})\cdot\ldots
\cdot(1-\ell^{-2g})}=\frac12.
$$
\end{proof}

\begin{rem}
Note that the factor $1/2$ fully agrees with the proof of Proposition \ref{agg}.
\end{rem}

We have described the geometry of $\b_1 \setminus \b_3$ in Section
\ref{sec:geomb3}. If we work with the level cover $\Perfl$, then the
irreducible components of $\b_1(\ell)$ are compactifications of the
universal level families $\calX_{g-1}(\ell)$, the components of the
locus $\Delta(\ell)$ are
${\calX_{g-2}(\ell)}\times_{\calA_{g-2}(\ell)}
\hat{\calX}_{g-2}(\ell)$, and $\b_2(\ell) \setminus \b_3(\ell)$ is
obtained by the gluing of $\ell$ copies of the $\ell$'th powers of
the universal Poincar\'e bundle (see \cite{gl} for a discussion of
why this is $\calP^{\otimes\ell}$, but this does not matter for our
purposes here). With this description of the geometry, we can prove
\begin{lm}
The second term in formula (\ref{levelintersection}) can be reduced
to a computation on $Y$ as follows: for fixed $a, b>0$ with $a+b=N$
$$
  \frac{\ell^N}{\nu_g(\ell)}\left\langle \sigma^*L^{G-N}
 \sum\limits_{i > j}D_i^a D_j^b \right\rangle_\Perfl
$$
$$
 =\frac18\left\langle L^{G-N}(-2\xi -2f^*T_1+f^*\calP)^{a-1}(-2\xi -2f^*T_2+f^*\calP)^{b-1}
\right\rangle_{Y}.
$$
\end{lm}
\begin{proof}
We denote by $d_{II}(\ell)$ the number of components of
$\mathop{\cup}\limits_{i>j} (D_i\cap D_j)$ --- note that all these
components lie in $\b_2(\ell)\subset\Perfl$ and that each component
of this set lies on only two boundary divisors, i.e. it is impossible
to have $D_i\cap D_j\subset D_k$ for $i,j,k$ distinct. This follows
since the second Voronoi decomposition for $g=2$ is basic. For any
fixed $i$ we think of $D_i$ as being a compactification $\XPl$ of
$\Xl$. The intersections $D_i\cap D_j$, being the closures of torus
orbits in a toroidal compactification, are irreducible, and are the
boundary components of $\XPl$. The fibers of the map $\Xpartl \to
\calA_{g-1}^{\op{part}}(\ell)$ are corank $1$ degenerations of
abelian varieties and as such have $\ell$ irreducible components.
Hence the number of irreducible boundary components of $\Xpartl$ is
equal to the number of irreducible boundary components of $\Perfl$
multiplied by $\ell$. Thus
$$
  d_{II}(\ell)=\frac{\ell}{2}d_g(\ell)d_{g-1}(\ell)
$$
where the factor $1/2$ accounts for the fact that every intersection
$D_i \cap D_j$ is contained in exactly $2$ boundary components. We
also need to compute $e_{II}(\ell)$, the degree of the map from any
$D_i\cap D_j$ to $\b_2$. Indeed, from our description of the geometry
we see that the level cover over $\b_2\setminus\b_3$ is the
composition of mapping $\ell$ copies of the $\PP^1$-bundle to one
copy (so factor of $\ell$ for the degree) together with mapping each
component of $\Delta(\ell)$ to $\Delta$, i.e. mapping
${\calX_{g-2}(\ell)}\times_{\calA_{g-2}(\ell)}
\hat{\calX}_{g-2}(\ell)\to {\calX_{g-2}}\times_{\calA_{g-2}}
\hat{\calX}_{g-2}$. The latter map has the degree equal to
$\nu_{g-2}(\ell)$, which is the degree of the map of the base
$\calA_{g-2}(\ell)\to\calA_{g-2}$, times the degree on the fibers,
i.e. the square of the degree of the map of a $\ell$-times
principally polarized abelian variety to a principally polarized
abelian variety. Thus we get
$$
  e_{II}(\ell)=\ell\nu_{g-2}(\ell)\ell^{4(g-2)}.
$$
Similarly to the first term above, to compare the intersection
numbers on $\Perfl$ and on $\Perf$ we need to compute the ratio
$$
 \frac{\ell^2d_{(II)}(\ell)e_{(II)}(\ell)}{\nu_g(\ell)}=
 \frac{\ell^2\ell d_g(\ell)d_{g-1}(\ell)\ell\nu_{g-2}(\ell)\ell^{4(g-2)}}
 {2\nu_g{(\ell)}}
$$
$$
  =\frac{\ell^{4+2g+2(g-1)+(g-2)(2g-3)+4(g-2)}(1-\ell^{-2g})(1-\ell^{2-2g})
  (1-\ell^{-2})\cdot\ldots\cdot(1-\ell^{4-2g})}{8\ell^{g(2g+1)}
  (1-\ell^{-2})\cdot\ldots\cdot(1-\ell^{-2g})}
$$
$$
=\frac18.
$$
This accounts for the factor $1/8$ which occurs on the right hand
side of the equation. The claim now follows from Lemma
\ref{lm:symmetricintersection}.
\end{proof}
The third term in (\ref{levelintersection}) can also be dealt with in
a similar fashion.
\begin{lm}\label{DDDisdelta}
Assume that $N < 3g-3$ (in which case $\b_3$
can be ignored). Then for fixed $a, b, c >0$ with $a+b+c=N$ we have
$$
  \frac{\ell^N}{\nu_g(\ell)}\langle \sigma^*L^{G-N}
 \sum\limits_{i > j >k}D_i^a D_j^bD_k^c \rangle_\Perfl
$$
$$
 =\frac{1}{12} \left\langle L^{G-N}(-2T_1 - \calP)^{a-1}
(-2T_2 - \calP)^{b-1}\calP^{c-1} \right\rangle_\Delta.
$$
\end{lm}
\begin{proof}
The intersection of three boundary components  $D_i \cap D_i \cap
D_k$ for different indices $i$, $j$ and $k$ can either be of local or
of global type. In the latter case such an intersection lies entirely
in $\b_3(\ell)$ and is irrelevant for us. Any intersection $D_i \cap D_i
\cap D_k$ of local type corresponds to a $2$-dimensional cone in the
perfect cone decomposition for $g=2$ and is thus a torus orbit and
hence irreducible. Moreover the intersections $D_i \cap D_i \cap D_k$
of local type are exactly the boundary components of $\Delta(\ell)$
and no component of $\Delta(\ell)$ is contained in $4$ boundary
components by Lemma \ref{no4intersect}. Similarly to the above case,
we denote by $d_{III}(\ell)$ the number of irreducible components of
$\mathop{\cup}\limits_{i> j > k} D_i\cap D_j\cap D_k$, where the
union is only taken over the intersections of local type. This is
equal to the number of irreducible components of $\Delta(\ell)$. For
fixed $i$, the number of components of $\Delta(\ell)$ contained in
$D_i$ is equal to the number of boundary components of ${\XPl}$ times
$\ell$ --- this is as in the second term computation above, since the
components of $\Delta(\ell)$ correspond exactly to the intersection
of two boundary components of  ${\XPl}$ and any such component
contains $2$ components of $\Delta(\ell)$. We thus have
$$
 d_{III}(\ell)=\frac{\ell}{3}d_g(\ell)d_{g-1}(\ell)
$$
where the denominator $3$ comes from the fact that every component of
$\Delta(\ell)$ lies in exactly $3$ boundary components. We also need
to compute $e_{III}(\ell)$, the degree of the map of each component
of $\Delta(\ell)$ to $\Delta$. The degree $e_{II}(\ell)$ that we
computed above was for mapping $\ell$ copies of a $\PP^1$-bundle to
one copy, together with some map of the base. Now we have the same
map of the base, and are mapping just one point taken from the union
of $\ell$ copies of $\PP^1$ in the fiber to one point, and thus
$e_{III}(\ell)=e_{II}(\ell)/\ell$. We compute
$$
 \frac{\ell^3d_{III}(\ell)e_{III}(\ell)}{\nu_g(\ell)}=
 \frac{\ell^3\ell d_g(\ell)d_{g-1}(\ell)\nu_{g-2}(\ell)\ell^{4(g-2)}}
 {3\nu_g(\ell)}=\frac{1}{12}
$$
since this is the same computation as we had for the second term, but
with 3 instead of 2 in the denominator. The claim now follows from
Lemma \ref{lm:intthreeboundcomp}.
\end{proof}
Combining these results, we finally get the following
\begin{prop}\label{threeterms}
For $2g-1\le N<3g-3$ the intersection numbers are
$$
  a_N^{(g)}=\langle L^{G-N}D^N\rangle_\Perf
  =\frac12\langle L^{G-N}(-2\T)^{N-1}\rangle_\Xpart
$$
$$
  + \frac 18\sum\limits_{a+b=N, a,b>0} \binom{N}{a}
\langle L^{G-N}(-2\xi -2f^*T_2+f^*\calP)^{a-1}(-2\xi -2f^*T_1+f^*\calP)^{b-1} \rangle_{Y}
$$
$$
  +\frac{1}{12} \sum\limits_{a+b+c=N, a,b,c>0}
  \binom{N}{a,\ b,\ c}\left\langle L^{G-N}(-2T_2 - \calP)^{a-1}
  (-2T_1 - \calP)^{b-1}\calP^{c-1} \right\rangle_\Delta
$$
$$
  =:(I)+(II)+(III)
$$
where we have numbered the terms for future reference.
\end{prop}
\begin{rem}
The factors $1/2$, $1/8$ and $1/12$ are due to the the stackiness of
$\Perf$ and the toroidal geometry. The first factor comes from the
fact that the map of stacks $\X \to D$ has degree $2$. Similarly, the
factor $1/8$ comes from the fact that each component of $Y(l)$ lies
on two boundary components (accounting for one factor $2$) and each
of these boundary components accounts for a further factor $2$.
Finally, each components of $\Delta(\ell)$ is the intersection of $3$
boundary components of $\Perfl$. {}From this one would naively expect
a factor $1/3\cdot2^3=1/24$. Note, however, that $\X$ and $Y$ are
families of abelian varieties, resp. degenerate abelian varieties
where the fibers have an involution, whereas the situation is
different for $\Delta$, which is actually a substack of $\Perf$
unlike $\X$ or $Y$. Hence we obtain the factor $1/12$ rather than
$1/24$. However, from these considerations it is not immediately
clear to us how to see that we have described the entire stabilizers
of these loci, and thus we gave the rigorous proofs above using
the combinatorics of level covers.
\end{rem}
Thus we have expressed the intersection number on $\Perf$ as a sum of
3 intersections numbers on different loci (without level cover!).
Each of the three summands is the intersection of classes for which
we have an explicit geometric description. In the following sections
we will compute these three intersection numbers.

\section{Intersection numbers on $\Delta$ (term (III))}
\label{sec:Intnumbersdelta}
We are now ready to compute the intersection numbers appearing in
formula (\ref{levelintersection}) where we will start with term
(III). {}This requires the computation of certain intersection
numbers on $\Delta$, and we will now determine completely the
intersection theory of divisors on $\Delta$. We will in fact compute
all the intersection numbers
$$
  \langle L^{G-N} T_1^lT_2^m\calP^n \rangle_{\Delta}=
  \langle L^{G-N}h_*(T_1^lT_2^m\calP^n)\rangle_{\calA_{g-2}}
$$
where the equality follows from the projection formula
\cite[Proposition 8.3 (c)]{fu}. In the following sections we will
reduce the computations of the two other terms to working on
$\Delta$.
\begin{thm}\label{pushTP}
The pushforward $h_*(T_1^lT_2^m\calP^n)=0$ unless $l=m$ and
$l+m+n=2g-4$, in which case it is
\begin{equation}\label{TTP}
h_*(T_1^{g-k-2}T_2^{g-k-2}\calP^{2k})=(-1)^k\frac{(g-2)!
(2k)!(g-2-k)!}{k!}[\calA_{g-2}].
\end{equation}
\end{thm}
We start by computing the pushforwards of the intersections of the
divisors $T_1$ and $T_2$.
\begin{lm}\label{TTlemma}
The pushforward $h_*(T_1^lT_2^m)$ is zero unless $l=m=g-2$, in which
case it is equal to $((g-2)!)^2[\calA_{g-2}]$.
\end{lm}
\begin{proof}
The projection map $h$ in diagram (\ref{project}) decomposes as a
composition of $p_1$ and $pr_2$, resp. of $p_2$ and $pr_1$; moreover,
the divisors $T_1$ and $T_2$ are the pullbacks of the universal
$(g-2)$-dimensional theta divisor $\T_{g-2}\subset\calX_{g-2}$
(symmetric and trivialized along the zero section, as always) under
$pr_1$ and $pr_2$ respectively. Thus
$$
 h_*(T_1^l T_2^m)=(p_1\circ pr_2)_*(pr_1^*\T_{g-2}^l
 \cdot pr_2^*\T_{g-2}^m)=p_{1*}(pr_{2*}(pr_1^* \T_{g-2}^l)\cdot \T_{g-2}^m )
$$
again by the projection formula. Since flat pullback commutes with
proper pushforward \cite[Proposition 1.7]{fu} we have
$$
 pr_{2*}\circ pr_1^*=p_1^*\circ p_{2*}.
$$
Hence we can rewrite the pushforward above as
$$
  p_{1*}(pr_{2*}(pr_1^* \T_{g-2}^l)\cdot \T_{g-2}^m )=
  p_{1*}(p_1^*(p_{2*} \T_{g-2}^l) \cdot \T_{g-2}^m).
$$
The advantage of doing this is that the pushforward under $p_{2*}$,
is simply the pushforward of the universal symmetric theta divisor
trivialized along the zero section on the universal family
$\calX_{g-2}\to\calA_{g-2}$, which we know by (\ref{push}) to be zero
for $l>g-2$. On the other hand  $p_{2*} \T_{g-2}^l=0$ for $l < g-2$
for dimension reasons. Since this is symmetric in $l$ and $m$ the
only remaining possibility is $l=m=g-2$ in which case we get
$$
  p_{1*}(p_1^*(p_{2*} \T_{g-2}^{g-2}) \cdot \T_{g-2}^{g-2})
  =(g-2)! p_{1*}(p_1^*([\calA_{g-2}]) \cdot \T_{g-2}^{g-2} )
$$
$$=(g-2)! p_{1*}(\T_{g-2}^{g-2})
  =((g-2)!)^2 [\calA_{g-2}].
$$
\end{proof}
One could now proceed to compute directly the pushforwards of
intersections of $T_1$ and $T_2$ with $\calP$, but there is an easier
way. We will use the automorphism group of $\Delta$. In Section
\ref{sec:inttheorydeltaY} we have already introduced the shift
operator $s: \Delta \to \Delta$ defined by $s(z,b)=(z+b,b)$. Let
$$
  V_N:=(s^N)^*(T_1)=T_1 + N^2T_2 + N\calP
$$
where the last equality follows from repeated application of Lemma
\ref{lm:shiftoperator}.

\begin{lm}\label{lm:Vn}
For any $N,l,m$ the pushforward $h_*(V_N^lT_2^m)=h_*(T_1^lT_2^m)$.
\end{lm}
\begin{proof}
By Lemma \ref{lm:shiftoperator} we have
$V_N^lT_2^m=(s^N)^*(T_1^lT_2^m)$, and the result then follows since
$s^N$ is an automorphism of $\Delta$ satisfying $h\circ s^N=h$.
\end{proof}
We can now use this to complete the computation of all pushforwards
on $\Delta$.
\begin{proof}[Proof of theorem \ref{pushTP}]
We will perform induction in $n$ --- the power of the Poincar\'e
bundle in $h_*(T_1^lT_2^m\calP^n)$. First note that for $n=0$ we know
the result to be true --- this is the content of Lemma \ref{TTlemma}.
Now assume the result to be true if the power of the Poincar\'e
bundle is strictly less than $n$, and let us prove the result for $n$
and arbitrary powers $l$ and $m$. Using Lemma \ref{lm:Vn}, we see
that for all $A,B,N$
$$
 h_*((T_1+N^2T_2+N\calP)^A T_2^B)=\delta_{A,\, g-2}\delta_{B,\, g-2} ((g-2)!)^2
$$
where $\delta$ denotes Kronecker's delta. The pushforward on the
left-hand side is a polynomial in $N$, which thus has to be
constant. In particular, if we extract the $N^n$ term, it must be
identically zero, i.e. for all $n>0$ we have
\begin{equation}\label{ncterm}
  \sum\limits_{i=0}^{\lfloor n/2\rfloor} \binom{A}{i,\ n-2i,\ A-n+i}
  h_*(T_1^{A-n+i}T_2^{B+i}\calP^{n-2i})=0.
\end{equation}
Observe that this pushforward contains one term which we would like
to compute, namely that where $\calP^n$ appears, while all the other
terms have lower powers of $\calP$. In particular we know inductively
that unless $(A-n+i)+(B+i)+(n-2i)=A+B=2g-4$, and $A-n+i=B+i$ for some
(and thus for all) $i$, all the terms with lower powers of $\calP$
are zero, and thus (of course the binomial coefficient is non-zero)
the term with $h_*(T_1^{A-n}T_2^B\calP^n)$ is also zero unless
$A+B=2g-4$ and $A-n=B$. Thus we have inductively proven that the
zeroes claimed for terms with $\calP^n$ are there, and it remains to
obtain the formula for the only non-zero pushforward with $\calP^n$
(for the case of $n$ even --- for $n$ odd still all the pushforwards
are zero).

In this case we have $A-n=B$ and $2B=2g-4-n$, so denoting $n=2k$ we
get $B=g-2-k$, and $A=g-2+k$. Since the above equation determines
$h_*(T_1^{g-2-k}T_2^{g-2-k} \calP^{2k})$ uniquely, to prove that
this pushforward is equal to what we want, it is enough to verify
that plugging in the claimed formula for it yields zero for the sum
above. Indeed, by plugging in the values of the pushforwards in
(\ref{ncterm}) above we get
$$
\begin{aligned}
 &\quad \sum\limits_{i=0}^k \binom{g-2+k}{i,\ 2k-2i,\ g-2-k+i}
 (-1)^{k-i}\frac{(g-2)!(2k-2i)!(g-2-k+i)!}{(k-i)!}\\
 &=(g-2)!(g-2+k)!\sum\limits_{i=0}^k (-1)^{k-i}
 \frac{(2k-2i)!(g-2-k+i)!}{i!(2k-2i)!(g-2-k+i)!(k-i)!}\\
 &=(g-2)!(g-2+k)!\sum\limits_{i=0}^k (-1)^{k-i}
 \frac{1}{i!(k-i)!}=0
\end{aligned}
$$
where we have used the standard identity
$$
  \sum\limits_{i=0}^k (-1)^{k-i}\binom{k}{i}=0
$$
for the last equality.
\end{proof}
In order to compute term (III) we must compute intersection numbers
of the form
$$
\langle L^{G-N}(-2 T_2 -\calP)^{a-1}(-2T_1 - \calP)^{b-1} \calP^{c-1} \rangle_{\Delta}
$$
$$
=\langle L^{G-N}h_*((-2T_2 - \calP)^{a-1}(-2T_1 - \calP)^{b-1}\calP^{c-1})\rangle_{\calA_{g-2}}
$$
with $a+b+c=N$. For $a+b \leq 2g-1$ we set
$$
C_g^{a,b}:=(-1)^{a+b+g}(g-2)!\sum\limits_{i=0}^{\min(a-1,b-1)}
 \frac{(-4)^i(a-1)!(b-1)!(2g-4-2i)!}{i!(a-1-i)!(b-1-i)!(g-2-i)!}
$$
and $C_g^{a,b}:=0$ otherwise.

This sum, as a generalized hypergeometric sum, can be evaluated
explicitly, using Maple \cite{maple}, which implements the method
described for example in \cite{concrmath}, to yield the following
answer
\begin{equation} \label{equ:C_g^{a,b}}
 C_g^{a,b}=(-1)^{a+b+g}(2g-4)!\frac{\Gamma\left(\frac52-g\right)
 \Gamma\left(\frac12+a+b-g\right)}{\Gamma\left(\frac32+a-g\right)
 \Gamma\left(\frac32+b-g\right)}
\end{equation}
where by using the standard properties of the $\Gamma$ function one can
further substitute
\begin{equation}\label{Gamma}
 \Gamma\left(\frac12+n\right)=\frac{(2n-1)!!}{2^n}\sqrt\pi\qquad
 {\rm and}\qquad \Gamma\left(\frac12-n\right)=\frac{(-2)^n}{(2n-1)!!}\sqrt\pi
\end{equation}
for $n\in\ZZ_{\ge 0}$.

With this notation we have the following
\begin{cor}
\label{cor:pushforwardthreeterms} The pushforward
$$
h_*((-2T_2 - \calP)^{a-1}(-2T_1 - \calP)^{b-1}\calP^{c-1})=0
$$
unless $a+b+c=2g-1$, in which case it is
$$
  h_*((-2T_2 - \calP)^{a-1}(-2T_1 - \calP)^{b-1}\calP^{c-1})=C_g^{a,b}[\calA_{g-2}].
$$
\end{cor}
\begin{proof}
All the monomials appearing in these pushforwards are zero unless
$(a-1)+(b-1)+(c-1)=2g-4$, so the vanishing is easy, and we are left
to compute for $a+b+c=2g-1$, using formula (\ref{TTP})
$$
h_*((-2T_2 - \calP)^{a-1}(-2T_1 - \calP)^{b-1}\calP^{c-1})
$$
$$
 = (-1)^{a+b-2}\!\!\!\sum\limits_{i=0}^{\min(a-1,b-1)}
 \binom{a-1}{i} \binom{b-1}{i}
h_*((2T_2)^i(2T_1)^i\calP^{2g-4-2i})
$$
$$
 =(-1)^{a+b}\sum\limits_{i=0}^{\min(a-1,b-1)}\frac{(-1)^{g-i}4^i(a-1)!(b-1)!(g-2)!(2g-4-2i)!}
{i!(a-1-i)!(b-1-i)!(g-2-i)!}.
$$
\end{proof}
Combining this with the previous results we finally obtain
\begin{thm} \label{theo:formulaforIII}
The term (III) in (\ref{levelintersection}) is equal to zero in all
cases except $N=2g-1$ when it is
$$
 (III)=\frac{(-1)^{g+1}(4g^2-8g+7)(2g-4)!}{12(2g-1)}
 +\frac{16^g(g-2)!(g-1)!}{192(2g-1)}
$$
$$
 +\frac{(2g-1)!(2g-4)!}{12}\sum\limits_{a=1}^{2g-3}\sum\limits_{k=1}^{2g-2-a}
 \frac{(-1)^{g+a+k+1}\left[\frac32-g+a\right]_k}{(k+1)(2g-a-k-2)!a!\left[\frac52-g\right]_k}
$$
where $[z]_k:=z\cdot(z+1)\cdot\ldots\cdot(z+k-1)$ denotes the
so-called Pochhammer symbol.
\end{thm}
\begin{proof}
This is a straightforward calculation. Indeed, by the above corollary
all the pushforwards in
$$
 (III)=\frac{1}{12}  \sum\limits_{a+b+c=N, a,b,c>0} \binom{N}{a,\ b,\ c}
$$
$$\times\, \left\langle L^{G-N}h_*\left((-2T_2 - \calP)^{a-1}
  (-2T_1 - \calP)^{b-1}\calP^{c-1}\right) \right\rangle_{\calA_{g-2}}
$$
are zero unless $N=a+b+c=2g-1$, in which case we get
$$
  (III)=\frac{1}{12}\langle L^{(g-2)(g-1)/2}\rangle_{\calA_{g-2}}\sum\limits_{a+b+c=2g-1, a,b,c>0}
\frac{(2g-1)!}{a!b!c!}C_g^{a,b}.
$$
Using Maple to sum this (and simplifying by hand each time in between
summations) yields the claimed formula.
\end{proof}
\begin{rem}
The resulting expression can
also be rewritten in terms of the hypergeometric functions ${}_3F_2$.
\end{rem}

\section{Intersection numbers on $Y$ (term (II))}
For term (II) in  (\ref{levelintersection}) we must compute
intersection numbers of the form
$$
  \langle L^{G-N}(-2\xi -2f^*T_1 + f^*\calP)^a(-2\xi -2f^*T_2 + f^*\calP)^b  \rangle_{Y}.
$$
We will compute the intersection number on $Y$ by pushing the
computation down to $\Delta$ and then to $\calA_{g-2}$, using the
results of the previous section. Recall from diagram (\ref{project})
the composite map $\pi=h\circ f:Y\to\calA_{g-2}$. We will need to
compute the pushforwards under this map of various products of
divisors. In doing so, we shall distinguish between the cases when
$\xi$ is one of the factors or is not. Note that by (\ref{xisquare})
we can always assume the power of $\xi$ to be at most one. We first
make the following observation
\begin{lm}\label{lm:pifh}
For any class $x\in \CH^*(\Delta)$ the following identities hold:
$$
  \pi_*f^*x=h_*(f_*f^*x)=0 \quad{\rm and}\quad \pi_*(\xi f^*x)=
  h_*(x)\quad\in \CH^*(\calA_{g-2}).
$$
\end{lm}
\begin{proof}
Note that the fiber dimension of the map $f$ is 1. Hence the first
claim follows for dimension reasons. The second claim follows since
$\xi$ is a section of $f$.
\end{proof}

\begin{prop} \label{prop:pushxiTP}
The following equality holds:
$$
  \pi_*((-2\xi -2f^*T_1 + f^*\calP)^{a-1}(-2\xi -2f^*T_2 + f^*\calP)^{b-1})\hskip5cm
$$
$$
  \hskip5cm=
  \begin{cases}
  2C_g^{a,2g-1-a}[\calA_{g-2}]& \mbox{if $a+b=2g-1$}\\
  0 & \mbox{otherwise.}
  \end{cases}
$$
\end{prop}
\begin{proof}
We compute
$$
  (-1)^{a+b}\pi_*\left((-2\xi -2f^*T_1 + f^*\calP)^{a-1}(-2\xi -2f^*T_2 + f^*\calP)^{b-1}\right)
$$
$$
  =\pi_*\left(\sum\limits_{i=0}^{a-1}
\sum\limits_{ j=0}^{b-1}\binom{a-1}{i}\binom{b-1}{j}(2f^* T_1)^i (2f^* T_2)^j(2\xi-f^*\calP)^{a+b-2-i-j}\right).
$$
Denoting $A:=a+b-2-i-j$ and using (\ref{xisquare}) repeatedly we
obtain
$$
 (2\xi-f^*\calP)^A=(\xi+(\xi-f^*\calP))^A
 =\xi^A+(\xi-f^*\calP)^A
$$
$$
 =\xi^A+(\xi-f^*\calP)^{A-1}(\xi-f^*\calP)
 =\xi^A+(\xi-f^*\calP)^{A-1}(-f^*\calP)=\ldots
$$
$$
 =\xi^A+(\xi-f^*\calP)(-f^*\calP)^{A-1}
 =(1+(-1)^{A-1})\xi(f^*\calP)^{A-1}-(f^*\calP)^A.
$$
We now apply Lemma \ref{lm:pifh} to get
$$
  \pi_*((2f^* T_1)^i (2f^* T_2)^j(2\xi-f^*\calP)^A)
$$
$$
  =\pi_*\left(\left((1+(-1)^{A-1})\xi(f^*\calP)^{A-1}
  -(f^*\calP)^A\right)(2f^*T_1)^i (2f^* T_2)^j\right)
$$
$$
 =h_*((1+(-1)^{A-1})\calP^{A-1}(2T_1)^i (2T_2)^j).
$$
By Theorem (\ref{pushTP}) this pushforward is non-zero only if $i=j$
and $a+b=2g-1$, in which case it is equal to
$$
  h_*((1+(-1)^{2g-4-2i})\calP^{2g-4-2i}(2T_1)^i (2T_2)^i)
$$
$$
 =2(-1)^{g-2-i}4^i\frac{(g-2)!(2g-4-2i)!i!}{(g-2-i)!}[\calA_{g-2}].
$$
Substituting this into the sum above gives the desired result.
\end{proof}
We can now finish the computation of the second term for the
intersection number --- this is just a matter of combining the known
results, and computing carefully.
\begin{thm}\label{theo:formulaforII}
Term (II) in (\ref{levelintersection}) is equal to
$$
(II)=\begin{cases}
\frac{-a_0^{(g-2)}}{64(2g-1)} \left(2^{4g}(g-1)!(g-2)!+32(-1)^g(2g-3)!\right)&\mbox{ if $N=2g-1$}\\
0 &\mbox{otherwise.}
\end{cases}
$$
where the explicit value of $a_0^{(g-2)}$ can of course be
substituted from formula (\ref{int0}).
\end{thm}
\begin{proof}
By definition, (II) is a sum involving the terms $\pi_*((-2\xi
-2f^*T_2+f^*\calP)^{a-1}(-2\xi -2f^*T_1+f^*\calP)^{b-1})$. Since
these are all zero unless $a+b=2g-1$ we obtain the vanishing of term
(II) as claimed. It follows from the definition of term (II) and
Proposition \ref{prop:pushxiTP} that for $N=2g-1$
$$
 (II)= \frac18 \langle L^{(g-2)(g-1)/2}\rangle_{\calA_{g-2}}
 \sum\limits_{ a+b=2g-1,a,b>0}\frac{(2g-1)!}{a!b!}\cdot
 2C_g^{a,b}.
$$
Note now that $C_g^{a,b}$ is symmetric in $a$ and $b$. Thus it is
enough to compute the sum for $a\le b=2g-1-a$, i.e. we have
$$
 (II) =\frac{1}{2} \sum\limits_{a=1}^{g-1}\binom{2g-1}{a}C_g^{a,2g-1-a}a_0^{(g-2)}.
$$
We now use Maple to sum the terms given explicitly by formula
(\ref{equ:C_g^{a,b}}), using the explicit factorial expressions for
$\Gamma$-functions (for each $\Gamma$-function factor the sign of $n$
in (\ref{Gamma}) is the same for all terms in the sum). This
summation yields the result as claimed.
\end{proof}

\section{Intersection numbers on $\Xpart$ (term I)}
In this section we finish the computation of the number
$a_{2g-1}^{(g)}$, together with the proof that the next $g-2$ numbers
after it are zero, by computing term (I) in formula
(\ref{levelintersection}). As in the previous section, we eventually
reduce this computation to a computation on $\Delta$ and thus get the
resulting number. The crucial ingredient in this reduction is a use
of the Grothendieck-Riemann-Roch theorem. Recalling that in the range
$N < 3g-3$ which we are interested in, $L^{G-N}$ is zero on $\b_3$,
and that $L$ is a pullback from the Satake compactification, we get
from Lemma \ref{lm:termIred}
\begin{equation}\label{Iispush}
(I)=\frac12(-2)^{N-1}\langle L^{G-N}\T^{N-1}\rangle_\Xpart
=\frac12(-2)^{N-1}\langle L^{G-N}\pi_*(\T^{N-1})\rangle_\Prt
\end{equation}
where $\pi:\Xpart\to\Prt$ is the partial compactification of the universal family in genus $g-1$.

Thus we need to compute the pushforward of $\T$ on the partial
compac\-ti\-fi\-ca\-tion of the universal family ---
note that formula (\ref{push}) only applies on the open part
$\X\to\A$. Unfortunately we do not know a way to
use a numerical argument similar to the one used in \cite{fc} and
explained in \cite{vdg2} to compute the pushforward on
$\Prt$ directly, and thus will use the following
identity
\begin{equation}\label{TFD}
 \pi_*(e^\T F)=e^{D/8}
\end{equation}
where $D$ denotes the class of the boundary of $\Prt$ and
$F:=\op{Td}^{\vee{}}(\calO_{{\Delta}})^{-1}$. Recall that ${\Delta}$
is the class on the stack $\Xpart$ which for all levels $n \geq 3$ is
given by the union of the singular loci of the fibers of the extended
universal family $\pi: \Xpart(n) \to \Prt(n)$. Formula (\ref{TFD})
was proven in \cite[Theorem 4.9]{vdg2} in cohomology, but using
\cite[Theorem 5.2]{ev} it follows that it also holds in the Chow
ring. This identity allows us to express term (I) in terms of the
pushforwards that we already know.

To simplify the notation, we denote the coefficients of the power
series expansion of the Todd class
$$
 \sum\limits_{i=0}^\infty b_n x^n:=\frac{x}{1-e^{-x}}=1+\frac{x}{2} +
\sum\limits_{k=1}^\infty(-1)^{k-1}\frac{B_k}{(2k)!}x^{2k}
$$
where $B_k$ are Bernoulli numbers, i.e. we set $b_0=1,
b_1=\frac{1}{2}$ and for $k\ge 1$ let  $b_{2k+1}=0$, and
$b_{2k}=(-1)^{k-1}\frac{B_k}{(2k)!}$.

\begin{lm}
Let $Z=E \cap F$ be a complete intersection of codimension $2$ on a
smooth manifold $X$. Then
\begin{equation}\label{horrible}
 \op{Td}^\vee((\calO_Z)^{-1})=\sum\limits_{i,j,k\ge 0}\sum\limits_{a=0}^k
\frac{(-1)^kb_ib_j}{(k+1)a!(k-a)!}E^{i+a}F^{j+k-a}.
\end{equation}
\end{lm}
\begin{proof}
We have an exact sequence
$$
 0 \longrightarrow \calO_X(-E-F) \longrightarrow \calO_X(-E)
 \oplus \calO_X(-F) \longrightarrow \calI_Z \longrightarrow 0.
$$
It then follows from the multiplicativity of the Todd class that
$$
 \op{Td}((\calO_Z)^{-1}) = \op{Td}(\calI_Z) =
\left(\frac{-E}{1-e^E}\right)\left(\frac{-F}{1-e^F}\right)
\left(\frac{-E-F}{1-e^{E+F}}\right)^{-1}
$$
and thus (note that taking the dual means changing the signs of $E$
and $F$)
\begin{equation}\label{formulatodd}
 \op{Td}^\vee((\calO_Z)^{-1})=\left(\frac{E}{1-e^{-E}}\right)\left(\frac{F}{1-e^{-F}}\right)
\left(\frac{1-e^{-E-F}}{E+F}\right).
\end{equation}
The claim then  follows by expanding the power series.
\end{proof}
\begin{rem}
Notice that this formula is written in terms of $E$ and $F$; however,
in each degree it is a symmetric polynomial of $E$ and $F$, and as
such is expressible in terms of the elementary symmetric polynomials
$Z=EF$ and $c_1(N_{Z/X})=E+F$. It is in this sense that we understand
this formula.
\end{rem}
Let us further denote by
$$
 b_{n,m}:=\sum\limits_{i=0}^n\sum\limits_{j=0}^m \frac{(-1)^{n+m-i-j}b_ib_j}{(n+m-i-j+1)(n-i)!(m-j)!}
$$
the coefficient of $E^nF^m$ in (\ref{horrible}). We were surprised to
find out that these coefficients admit a remarkably simple formula,
which we have been unable to find in the literature.
\begin{lm}\label{lm:bnm}
The coefficients above can be computed as follows:
$$
  b_{n,m}=\begin{cases}(-1)^{\frac{n-m}{2}}\frac{B_{\frac{n+m}{2}}}{(n+m)!}&
  \mbox{if $n+m$ is even and $nm\ne0$}\\ 1&\mbox{if $n=m=0$}\\ 0&\mbox{otherwise.}\end{cases}
$$
\end{lm}
\begin{proof}
We start by noting the following straightforward identity:
$$
 \frac{1-e^{-E-F}}{(1-e^{-E})(1-e^{-F})}=\frac{1}{1-e^{-E}}+ \frac{1}{1-e^{-F}} -1.
$$
Multiplying this by $EF$, we recognize the expressions for the Todd
classes of $E$ and $F$ (multiplied by $F$ and $E$, respectively) on
the right. Finally we divide by $E+F$, getting expression
(\ref{formulatodd}) on the left, and use
$$
 \frac{E^{2k}+F^{2k}}{E+F}=E^{2k-1}F-E^{2k-2}F^2+\ldots-E^2F^{2k-2}+EF^{2k-1}
$$
on the right (note that there are no linear terms here, unlike
the Todd class).
\end{proof}
We now want to apply this to $\Delta$ considered as a codimension 2
substack of $\Xpart$. This is a complete intersection in a stack
sense: if we go to level covers then every component of
$\Delta(\ell)$ in $\Xpartl$, which we can identify with a fixed
boundary component $D_i$, is a complete intersection of two further
boundary components $D_j$ and $D_k$. By Proposition
\ref{normaltoDelta} we have
$$
  N_{\Delta/\Xpart}= \calP \oplus (\calP^{-1} \otimes T_2^{-2})
$$
and thus we can use $E=\calP$ and $F=-2T_2-\calP$ in the corollary
above.

\begin{prop}
The first term in (\ref{threeterms}) is zero unless $N=2g-1$, in
which case it is equal to
$$
 (I)=\frac{(2g-2)!}{2^g(g-1)!}a_{g-1}^{(g-1)}+(-1)^g(2g-3)!a_0^{(g-2)}
 \sum\limits_{m=1}^{g-1}\frac{(-1)^m2^{2m+2-2g}B_m}{(2g-2m-1)!(2m)!}
$$
\end{prop}
\begin{proof}
Expanding identity (\ref{TFD}) in power series and extracting the
degree $N-1$ term on the left, we get
$$
 \frac{{\pi}_*(\T^{N-1})}{(N-1)!}+
 \sum\limits_{k=1}^{N-1} {\pi}_*
 \left(\frac{\T^{N-1-k}}{(N-1-k)!}
 \Delta\sum\limits_{n=1}^{k-1}b_{n,k-n}\calP^{n-1}(-2T_2-\calP)^{k-n-1}\right)
$$
$$
=\frac{D^{N-g}}{8^{N-g}(N-g)!}.
$$
Computing an intersection number of any divisor class with $\Delta$
is the same as restricting to the locus $\Delta\subset\Perf$. {}From
Corollary \ref{thetaondelta} we know that $\T|_\Delta=T_1 +\frac12
\calP$. The map ${\pi}$ restricted to $\Delta$ is the map
$pr_1:\Delta\to\calX_{g-2}=\Prt\setminus\A$, and $L$ on $\Delta$ is a
pullback from $\calA_{g-2}$, so we are in the setup of the previous
sections, and know all the pushforwards from Theorem (\ref{pushTP}).
We can thus compute
$$
 (I)=(-1)^{N-1}2^{N-2}\left\langle L^{G-N}{\pi}_*(\T^{N-1})\right\rangle_\Prt
$$
$$
=(-1)^{N-1}2^{N-2}(N-1)!\left[\left\langle L^{G-N}
 \frac{D^{N-g}}{8^{N-g}(N-g)!}\right\rangle_\Prt \right.
$$
$$
-\sum\limits_{k=1}^{N-1}
 \left\langle \frac{L^{G-N}}{2}\sum\limits_{n=1}^{k-1}
 \frac{b_{n,k-n}}{(N-1-k)!}\right.
$$
$$
\left.\left.
\times\, h_*\left((T_1+\calP/2)^{N-1-k}
 \calP^{n-1}(-2T_2-\calP)^{k-n-1}\right)\right\rangle_{\calA_{g-2}}\right].
$$
Here the factor $1/2$ in the expression $L^{G-N}/2$ comes from the
fact that we consider $\Delta$ as a stack, i.e. we have an extra
involution which acts trivially on the underlying variety (cf. also
the proof of \cite[Lemma 4.8]{vdg2}).

Since all the pushforwards on the right are zero unless $N=2g-1$, and
so is the first term on the right, it follows that (I) is zero unless
$N=2g-1$. In this case we deduce from Corollary
\ref{cor:pushforwardthreeterms}
$$
 h_*\left((T_1+\calP/2)^{2g-2-k}\calP^{n-1}(-2T_2-\calP)^{k-n-1}\right)
=\frac{(-1)^{-k-1}}{2^{2g-2-k}}C_g^{k-n,2g-k-1}
$$
and thus get in that case
$$
(I)=2^{2g-3}(2g-2)!\left(\frac{a_{g-1}^{(g-1)}}{8^{g-1}(g-1)!}- \right.
$$
$$
\left. \frac{1}{2}a_{0}^{(g-2)}
\sum\limits_{k=1}^{2g-2}\sum\limits_{n=1}^{k-1}
 \frac{(-1)^k b_{n,k-n}}{2^{2g-2-k}(2g-2-k)!}C_g^{k-n,2g-k-1}\right).
$$
We can now use the explicit expression for the coefficients $b_{n,m}$
obtained in Lemma \ref{lm:bnm}. It follows that the coefficient
$b_{n,k-n}$ in the formula above is zero unless $k=2m$, in which case
$b_{n,2m-n}=(-1)^{m+n}\frac{B_m}{(2m)!}$, so that the double sum
above simplifies to
$$
 \sum\limits_{m=1}^{g-1}\frac{(-1)^mB_m}{2^{2g-2-2m}(2g-2-2m)!(2m)!}\sum\limits_{n=1}^{2m-1}
(-1)^n C_g^{2m-n,2g-2m-1}.
$$
We then use Maple to compute the sum over $n$ explicitly, getting a
simple expression $-\frac{(2g-3)!}{2g-1-2m}$ for it (this is again an
application of the algorithm for computing sums with hypergeometric
terms described in \cite{concrmath}), and end up with the expression
as claimed.
\end{proof}
\begin{cor}\label{theo:formulaforI}
The first term can be expressed in closed form as follows:
\begin{equation}
 (I)=\frac{(2g-2)!}{2^g(g-1)!}a_{g-1}^{(g-1)}-(-1)^g2^{2-2g}(2g-3)!a_0^{(g-2)}
\end{equation}
$$
 \times\left({\rm the\ coefficient\ of}\ x^{2g-1}{\rm\ in}\
 \left(\frac{2x}{1-e^{-2x}} -1-x\right)\frac{e^x-1}{x}\right)
$$
\end{cor}
\begin{proof}
Indeed, we recall that by definition $\sum\limits_{k=0}^\infty
b_k(2x)^k=\frac{2x}{1-e^{-2x}}$. If we subtract from this power
series the linear part, i.e. subtract $1+x$, and multiply the result
by $\frac{e^x-1}{x}= \sum\limits_{k=0}^\infty \frac{x^k}{(k+1)!}$,
the coefficient of $x^{2g-1}$ in the product is exactly equal to
$2^{2-2g}$ times the sum in the right-hand-side of the proposition.
\end{proof}

Combining this with the computation of terms (II) and (III) in the
previous section, we see that the only non-zero $a_N^{(g)}$ in our
range is the one for $N=2g-1$, and the expression for
$a_{2g-1}^{(g)}$, being the sum of three terms (I), (II), and (III),
is as claimed.

\section{Numerical examples}
\label{sec:examples}
Here we list the numerical values of the terms (I),(II) and (III) as
well as $a_{2g-1}^{(g)}$ for small values of $g$ (obtained by
evaluating our formulas using Maple)
\begin{equation}
 \begin{array}{|c|r|r|r|r|}
 \hline
 {}_{\rm genus}\!\!\backslash^{\rm term}&\hfill (I)\hfill &\hfill (II)\hfill &\hfill (III)\hfill &\hfill a_{2g-1}^{(g)}\hfill\\ \hline
 \vphantom{\dfrac{1}{2}}2  & \frac{1}{12}  & -\frac{3}{2} & \frac{1}{2} & -\frac{11}{12}\\
 \vphantom{\dfrac{1}{2}}3  & -\frac{1}{80} & -\frac{25}{24} & \frac{5}{24} & -\frac{203}{240}\\
 \vphantom{\dfrac{1}{2}}4  & \frac{1}{672} & -\frac{49}{80} & \frac{7}{80} & -\frac{1759}{3360}\\
 \vphantom{\dfrac{1}{2}}5  & -\frac{1}{1296} & -\frac{3637}{2520} & \frac{1063}{7560} & -\frac{59123}{45360}\\
 \vphantom{\dfrac{1}{2}}6  & \frac{1}{220} & -\frac{23837}{315} & \frac{1639}{315} & -\frac{976649}{13860}\\
 \vphantom{\dfrac{1}{2}}7  &-\frac{11}{18}&-\frac{4194073}{189}&\frac{17594928013}{16329600}&-\frac{49254708341}{2332800}\\
 \hline
 \end{array}
\end{equation}
Note that the terms for $g=2,3$ and $4$ coincide with \cite[Table
2b]{vdg1}, \cite[Table 3d]{vdg1} and \cite[Theorem 1.1]{egh},
respectively. For the last number we have to take an extra factor of
$1/2$ into account for $g=4$ since the computation in \cite{egh} was
done for varieties rather than stacks. The values for $g>4$ are new.

\section{Further comments}
\label{sec:comments}
Throughout the paper we have worked with the
perfect cone compactification $\Perf$ of $\calA_g$. In this case the
Picard group has two generators, namely the Hodge line bundle $L$ and
the boundary $D$. In particular, the boundary is irreducible and a
Cartier divisor (on the stack). There are other toroidal
compactifications, such as the Voronoi compactification $\Vor$ and
the Igusa compactification $\Igu$. Their Picard groups are much more
complicated. Nevertheless, our computations are also relevant for
other toroidal compactifications.

In general, a toroidal compactification is determined by the choice
of an admissible fan $\Sigma(F)$ for each cusp $F$. These choices
must be made in such a way that the resulting collection of fans
$\tilde{\Sigma}= \{ \Sigma(F)\}$ fulfills certain compatibility
conditions. In the case of $\calA_g$ there exists exactly one cusp in
each dimension, corresponding to the strata of the Satake
compactification $\Sat$ (cf. (\ref{eq:sat})). This means that, due to
the compatibility conditions, a toroidal compactification of
${\cal{A}}_g$ is already determined by the choice of an admissible
fan in the cone ${\Sym}^{\geq 0}(g,\RR)$ of semi-positive definite
real symmetric $g \times g$-matrices. Given any toroidal
compactification $\Abar$,  the geometry of the partial
compactification $\Abar \setminus \b_k$ is determined by the
intersection of the fan $\Sigma$ with ${\Sym}^{\geq 0}(k-1,\RR)$
where we can think of elements in ${\Sym}^{\geq 0}(k-1,\RR)$ as $(g
\times g)$-matrices by putting the matrices in the left hand corner
and extending this by zeroes in the other entries. All admissible
fans in ${\Sym}^{\geq 0}(g,\RR)$ coincide for $g \leq 3$ (of course,
up to taking a subdivision of this fan, which will result in a
blow-up of the toroidal compactification, but we will disregard
this), in particular $\Abar \setminus \b_3= \Perf \setminus \b_3$. As
all our computations happen outside $\b_3$, they apply to all
toroidal compactifications (disregarding artificial blow-ups), where
the meaning of $D$ is that it is the closure of the boundary of
Mumford's partial compactification.

\medskip

\noindent
Cord Erdenberger, Klaus Hulek,\\
Institut f\"ur Algebraische Geometrie,\\
Leibniz Universit\"at Hannover\\
Welfengarten 1, 30060 Hannover, Germany\\
{\texttt erdenber@math.uni-hannover.de}, {\texttt
hulek@math.uni-hannover.de}\\

\medskip

\noindent
Samuel Grushevsky,\\
Mathematics Department,\\
Princeton University,\\
Fine Hall, Washington Road,\\
Princeton, NJ 08544, USA \\
{\texttt sam@math.princeton.edu}
\end{document}